\documentclass[12pt]{article}
\usepackage{amsmath}
\usepackage{amssymb}
\usepackage{mathptmx}
\usepackage{amsfonts}
\usepackage[mathscr]{eucal}
\usepackage[dvips]{graphicx}
\usepackage{color}

\newcommand{\mbf}[1]{\ensuremath{\mathbf{#1}}}
\newcommand{\ms}[1]{\ensuremath{\mathscr{#1}}}

\newcommand{\Hb}{H_{\textrm{b}}}
\newcommand{\Hpar}{H_{\textrm{p}}}
\newcommand{\HI}{H_{\textrm{I}}}

\newcommand{\tens}{\otimes}

\newcommand{\LtwoRd}{L^{2}(\mathbf{R}^{d})}
\newcommand{\LtwoRdx}{L^{2}(\mathbf{R}^{d}_{\mathbf{x}})}

\newcommand{\Rd}{\mathbf{R}^{d} }

\newcommand{\omegab}{ \omega_{\, \textrm{b}} }

\newcommand{\bos}{\textrm{b}}

\newcommand{\parti}{\textrm{p}}
\newcommand{\rel}{\textrm{rel}}
\newcommand{\real}{\textrm{real}}
\newcommand{\Euc}{{\small \textrm{E}}}

\newcommand{\Br}{\textrm{Br}}

\newtheorem{theorem}{Theorem}[section]
\newtheorem{proposition}[theorem]{Proposition}

\newtheorem{corollary}[theorem]{Corollary}
\newtheorem{remark}{Remark}[section]

\begin{document}
\begin{center}
{\LARGE Functional Integral Representation for Relativistic Schr\"{o}dinger Operator Coupled to a Scalar Bose Field
 with  $P(\phi)$ Interaction } \\
 $\;$ \\
 Toshimitsu Takaesu \\
 $\;$ \\
\textit{Faculty of Mathematics, Kyushu University,\\  Fukuoka, 812-8581, Japan }
\end{center}

\begin{quote}
\textbf{Abstract}. In this paper  the system of  a semi-relativistic particle interacting with  a scalar Bose field  is investigated. The ultraviolet cutoff condition is imposed on the Bose field.  In the main theorem, the functional integral representation of 
the semi group generated by the total Hamiltonian with $P(\phi)$ interaction is obtained.  \\
 
{ \small
 $\qquad \qquad \; \; \; \; $   L\'{e}vy process,  Gaussian random process. } 
\end{quote}

\section{Introduction}
In this paper   the  system of a quantum particle  interacting with a  scalar Bose field 
 is investigated. 
The partilcle's Hamiltonian  is given the 
 relativistic Schr\"{o}dinger operator  with potential
\begin{equation}
    H_{\parti} \;  = \;     \sqrt{-\triangle + M^2 }  \, - \,  M   \; + \; V \label{Hparti} 
\end{equation}
 on the Hilbert space $L^2 (\mbf{R}^d_{\mbf{x}})$, where $M>0 $ is the  rest mass of the particle. For the  stochastic analysis of $\Hpar $, 
 the asymptotic behavior of the eigenvector is analyzed by the 
 functional integral representation in \cite{CMS90}, and the functional representation in electromagnetic  potential is derived  in   \cite{HIL09, IcTa86}.  In this paper, we construct the functional integral representation according to \cite{HIL09}.  For other results on the stochastic analysis of $\Hpar$, refer to e.g. \cite{IcTa86, KaLo11, LM10}.

$\;$ \\ 
  A scalar Bose  field  is constructed by stochastic process which investigated in constructive quantum field theory
(refer to e.g. \cite{GJ87, Si74}). 
 The field operators $ \{ \phi (f)  \}_{f \in \ms{K_{\bos}} } $   are defined by the Gaussan random process indexed  by a Hilbert space $\ms{K}_{\bos}$ on a probability space   $(Q_{\bos}, \mathfrak{B}_{\bos} , P_{\bos} )$.
 The state space is given  by $    L^2 ( Q_{\ms{K}_{\bos}} ) $ and the free Bose Hamiltonian $H_{\bos} $ is defined 
  by the differential second quantization of $ \omegab (-i \nabla) $ where $\omegab $ is  non-negative and continuous function.   Physically $\omegab (\mbf{k}) \geq 0$ denotes the one-particle energy of the field with momentum $\mbf{k}$.
Thus the triplet $ \left( L^2 ( Q_{\ms{K}_{\bos}} )  \, , \Hb,   
 \{\phi (f) \}_{ f \in  \ms{K}_{\bos} } \right)$  of the scalar Bose field is  defined.

$\;$ \\
The system of semi-relativistic particles coupled to a scalar Bose field is defined as follows.  
The state space is given by  $\ms{H}  = \LtwoRdx \tens L^2 (Q_{\bos } ) \simeq \, \int^{\oplus}_{\mathbf{R}^d}  L^2 (Q_{\bos } ) d \mbf{x}$ where $ \int^{\oplus}$  denotes the fibre direct integral. The free Hamiltonian is defined by
$ H_{0} = \Hpar \tens I + I \tens \Hb $ and the total Hamiltonian by 
\begin{equation}
H_{\kappa} \; \; = \; \;  H_{0} \;  \dotplus   \; \kappa \int_{\Rd}^{\oplus} P (\phi ({\rho_{\mbf{x}}} ))  d \mbf{x}   
\end{equation}
where $\dotplus$ denotes the form sum,  $P(\lambda  ) = \sum\limits_{j=1}^{2n} c_{j} \lambda^j  $,
 $c_j \in \mbf{R}$, $j=1, \cdots , 2n-1$, $c_{2n} > 0$,  and the ultraviolet cutoff condition $ \rho_{x} \in \ms{S}_{\textrm{real}}'$ for each $\mbf{x} \in \Rd$ is supposed. 
  
$\;$ \\ 
By using the  functional integral representations of $e^{-t\Hpar}$ and $e^{-t \Hb}$,     the functional integral representation of $e^{-tH_{\kappa}}$ is derived in the main theorem. 
Then, from the functional integral representation, it is seen that $e^{-t H\kappa} $ is positivity improving. Then,  as a corollary of the main theorem,  it is seen that  the  ground state of $H_{\kappa} $ is  unique if it exists.
$\;$ \\

For the spectral analysis for quantum particles systems coupled to Bose fields by the methods of stochastic analysis has been analyzed.   
For  non-relativistic QED model, its functional representation is obtained in \cite{Hi97}, and the case  with spin is considered in \cite{HiLo08}.  
The self-adjointness of the Hamiltonian is investigated in \cite{Hi02}, the analysis of the bound state  in \cite{Hi00b} and that of the exponential decay  in \cite{HiHi10}.  
 For spin-boson model and the Nelson model, the applications of their functional integral representation to spectral analysis are investigated   in \cite{BHLMS02, BH09, Hi99,  HLT11, LMS02b, Sp89, Sp98}. 
 
$\;$ \\
This paper is organized as follows. In section 2, the functional integral representation for the semi-relativistic particles 
is constructed. In section 3, we overview the Euclidean quantum field theory, and the functional integral representation for scalar Bose field is derived. In section 4, the interaction system is introduced. The the main theorem  is stated, and its proof is given. 
\section{Relativistic Schr\"{o}dinger Operator}

According to \cite{HIL09}, the functional integral representation of semigroup  generated by the  relativistic schr\"{o}dinger operator is derived  as follows. In this derivation,  the L\'{e}vy subordinator plays an important role.  A stochastic process $\{T_t\}_{t\geq 0}$  is called a L\'{e}vy subordinator if $\{T_t\}_{t\geq 0}$  is one dimensional L\'{e}vy process  starting at zero and almost surely  non-decreasing in $t \geq 0$.  

The function $\Psi \in C^{\infty} ((0, \infty ) ) $  is called  Bernstein function if
$\Psi \geq 0 $ and $(-1)^n \frac{d^n \Psi}{d x^n } \leq 0 $ for all $n \in \mbf{N}$. 
It is known  in (\cite{HIL09} ; Proposition 2.5)  that for a  Bernstein function $\Psi $ satisfying 
$\lim\limits_{x \to +0} \Psi (x) = 0 $, there exists an unique L\'{e}vy subordinator
$  \{T_t^{\Psi} \}_{t\geq 0} $  such that $ \mathbb{E} [ e^{-sT_t^{\Psi}} ] \;  = \;  e^{-t \Psi (s )}  $.

$\;$ \\
Let $M > 0 $ be the fixed mass of the particle, and let us set $h_{\rel} ( s) \; = \; \sqrt{s + M^2} -M $, $s> 0$. Since $h_{\rel} $ is a Bernstein function, it is seen  that  there exists a L\'{e}vy subordinator $\{T_t \}_{t\geq 0}$ on the probability space
$(\Omega_{\rel} , \mathfrak{B}_{\rel} , P_{\rel} ) $ satisfying 
\begin{equation}
\mathbb{E}_{\rel} [ e^{-sT_t}] \;  = \;  e^{-t h_{\rel} (s )}   , \label{8/10.2.1}
\end{equation}
where $ \mathbb{E}_{\rel} [X] = \int_{\Omega_{ \rel}} X (\eta ) d P_{\rel} (\eta ). $
\begin{remark}
Let $Y_{s} = B_{s} + Ms $, $s>0 $ where  $ \{ B_{s} \}_{s>0}$ is one dimensional 
Brownian motion starting at zero. It is known that $ \{T_t\}_{t\geq 0 }$ is represented as the first hitting time process
$T_{t} = \frac{1}{2} \inf \{  s > 0 | Y_{s} = t \}$.  (See  \cite{LHB11}; Example 2.18).
\end{remark}

$\;$ \\
Let $\{ \mbf{B}_t \}_{t \geq 0 } $ be $d$-dimensional Brownian motion starting $\mbf{x} $ on the probability space
$( \Omega_{\Br} , \mathfrak{B}_{\Br} ,  P_{\Br}^{\mbf{x}} ) $.  We introduce the probability space
\[
 ( \Omega_{\parti} , \mathfrak{B}_{\parti} , P^{\mbf{x}}_{\parti} ) = 
 ( \Omega_{\rel}  \times \Omega_{\Br}  ,  \overline{\mathfrak{B}_{\rel} \times  \mathfrak{B}_{\Br}} ,
  \overline{ P_{\rel} \times P^{\mbf{x}}_{\Br}}  ) , 
\]
 and let us define a stochastic process $\{ \mbf{X}_t \}_{t \geq 0} $ on 
$( \Omega_{\parti}, \mathfrak{B}_{\parti}  , P_{\parti}^{\mbf{x}})$
 defined by 
 \[
 \mathbf{X}_{t} \left( 
\left[ \begin{array}{cc} \eta \\ \omega \end{array} \right] 
 \right)
 \; =  \; \mathbf{B}_{T_t (\eta )} (\omega )   .  
\]
For $\psi \in  \LtwoRdx $, let us set 
\[
U_t^{\parti} \psi (\mbf{x} )
\; = \;  \mathbb{E}_{\parti}^{\mbf{x}}   [ \psi ( \mbf{X}_t ) e^{- \int_0^t V (\mbf{X}_s ) ds   } ]  ,
\]
where $ \mathbb{E}^{\mbf{x}}_{\parti} [Z] = \int_{\Omega_{\parti}} Z (\xi ) dP_{\parti}^{\mbf{x}} (\xi ) $. Here we assume  the following condition. 
\begin{quote}
\textbf{(S.1)} $V\in L^{\infty} (\Rd ) $.
\end{quote}
By Fourier transform and (\ref{8/10.2.1}), it is seen that for  rapidly decreasing function $\psi \in \ms{S} (\Rd) $,
\[
(e^{-t  (\sqrt{ - \triangle + M^2} -M) } \psi) (\mbf{x}) = 
\int_{\Rd} e^{- t h_{\rel} (\mbf{k}^2) } \hat{\psi} (\mbf{k}) e^{i \mbf{k} \cdot \mbf{x} } d\mbf{k}
= \mathbb{E}_{\rel} [ \mathbb{E}_{\Br}^{\mbf{x}} [\psi (\mbf{B}_{T_t} )]] =
 \mathbb{E}_{\parti}^{\mbf{x} } [ \psi (\mbf{X}_t ) ]. 
\]
 Then from  Trotter-Kato product formula, the functional integral representation for the semi-relativistic particle  is obtained:

$\;$ \\
\textbf{Proposition A} (\cite{HIL09} ; Theorem 3.8)
\\ Assume \textbf{(S.1)}. Then  $ ( \phi , e^{-t \Hpar } \psi )  \; =   ( \phi, U_{t}^{p} \psi )   $ for $\phi, \psi \in \LtwoRd $.

$\;$ \\
As a remark, in the proof of (\cite{HIL09}; Theorem 3.8)  it is proven that a Feynman-Kac formula  
\begin{equation}
e^{-t_1 \Hpar} g_1e^{ -(t_2 -t_1 ) \Hpar} 
  \cdots g_{n-1} e^{ -(t_{n} -t_{n-1} ) \Hpar}  \psi (\mbf{x}) \; = \; 
  \mathbb{E}_{\parti}^{\mbf{x} } [   \prod_{j=1}^{n-1} g_j (\mbf{X}_{t_j} )  \psi (\mbf{X}_{t_n }) 
e^{-\int_0^{t_n }  V_{\parti} (\mbf{X}_s)  ds }  ] \label{FKS_relativistic}
\end{equation}
 holds where $g_{j}  \in  L^{\infty} (\Rd ) , j=1, \cdots , n -1 $ and $\psi \in \LtwoRd$.  \\

\section{Scalar Bose fields}
\subsection{Gaussian random process indexed by Hilbert space}
In this subsection, basic properties for Gaussian  random process are explained. 
To construct Bose fields, the following proposition is needed. (See, e.g.\cite{Arai10};Theorem 2.5, \cite{LHB11};Theorem 5.9  )

$\;$ \\
\textbf{Proposition B} (Existence of Gaussian random process) \\
Let $\ms{K}$ be a separable and real Hilbert space. Then there exist  a stochastic process $\{ X_f \}_{f \in \ms{K}}$  indexed by $ \ms{K}$ on a probability space $(Q_{\ms{K}}, \mathfrak{B}_\ms{K} , P_{\ms{K}} )$ satisfying the following conditions. \\
$\quad$(\textbf{G.1}) For all $f \in \ms{K} $, $X_f $ is Gaussian random variable satisfying
 $  \mathbb{E} [e^{-it \phi_{f}}]  = e^{-\frac{ \|  \, f\|^2}{4} t^2}$. \\
$\quad$(\textbf{G.2}) $X_{af + b g} = a X_f + b X_g  $  for all $ f, g \in \ms{K} $ and $ a, b \in  \mbf{R}$. \\
$\quad$(\textbf{G.3}) $\mathfrak{B}_{\ms{K}} $ is the minimal $\sigma$-field generated by 
$\{ X_f \}_{f \in \ms{K} } $.  \\
\begin{remark}
The stochastic process    $\{ X_f \}_{f \in \ms{K}}$  satisfying \textbf{(G.1)}-\textbf{(G.3)} is called 
 the Gaussian random process indexed by $\ms{K}$. 
\end{remark}
Let  $\{ X_f \}_{f \in \ms{K}}$   be  the Gaussian random process indexed by $\ms{K}$. 
Then it is seen from \textbf{(G.3)}, that $\ms{D}_{0, \ms{K} } = \left\{   F ( X_{f_1} , \cdots , X_{f_{n}}   ) 
 \left| \right.  F \in \ms{S}_{\textrm{real} } (\Rd ) , f_{j} \in \ms{K} , j= 1, \cdots, n , n \in \mbf{N}  \right\}$
is dense in $L^{2}{ (Q_{\ms{K}})}$. Let 
$ L^2_n ( Q_{ \ms{K} }   ) $  be the closure of the linear hull of  the  set $  
    \{   : \prod_{j=1}^n X_{f_j}  : \, \left. \right| f_j \in \ms{K}, \, j = 1, \cdots , n \}
\cup \{1  \} $ 
where $ : \prod_{j=1}^n X_{f_j}  :  $ denotes the wick product defined recursively by $ : \prod_{j=1}^n X_{f_j}  :  
= X_{f_1} : \prod_{j=2}^{n} X_{f_j} : - \frac{1}{2} \sum_{j=2}^{n}  (f_1 , f_j) : \prod_{j \ne l}   X_{f_{l}} : $ and $: X_{f}: = X_f $.  It is seen that $ L^2_j ( Q_{ \ms{K} }  ) \bot L^2_l ( Q_{ \ms{K} }  ) $ for $j \ne l $. It is known that the Winer-Ito-Segal decomposition  
$ L^2 ( Q_{ \ms{K} } ) =\bigoplus\limits_{n=0}^{\infty} L^2_n ( Q_{ \ms{K} }  ) $ follows (See e.g. \cite{Arai10}; Lemma 2.13 , \cite{LHB11}; Lemma 5.4). 
 Let  $S $ be a closed operator on $\ms{K}$.   $\Gamma (S ) = \oplus_{n=0}^\infty \Gamma^{(n)} (S ) $ is called the  second quantization of $S$ defined by   $ \Gamma^{(n)} (S )   X_{ f_{1}} \cdots   X_{f_{n}} \; = \;  : X_{S f_{1}} \cdots   X_{Sf_{n}}$ for $f_j \in \ms{D} (S)$, $j=1, \cdots , n $, $n  \geq 0$.   In addition,
$d \Gamma (S ) = \oplus_{n=0}^\infty d \Gamma^{(n)} (S )  $ is called the differential second quantization of $S$ defined by $d\Gamma^n  (S )   X_{f_{1}} \cdots   X_{f_{n}} \; = \;  \sum_{j=1}^n : X_{ f_{1}} \cdots    X_{S f_j }  \cdots X_{f_{n}} : $  where $f_j \in \ms{D} (S)$, $j=1, \cdots , n $,   $ n \geq 0$.

\subsection{Construction of a scalar Bose field}
Let 
\[
\ms{K}_{\bos}^0 \; = \; 
\left\{  f \in \ms{S}_{\real}' ( \Rd ) \left| \frac{}{} \right. \int_{\Rd} 
\frac{|\hat{f} (\mbf{k})|^2 }{\omega_{\, \bos} (\mbf{k})}  d \mbf{k} \; \infty  \right\} ,
\]
where   $\ms{S}_{\real}' ( \Rd )  $ is the  space of real-valued tempered distributions, and   set 
\begin{equation}
(g, f )_{\ms{K}_{\bos}} \; = \; 
 \int_{\Rd} 
\frac{ \overline{\hat{g} (\mbf{k})} \hat{f}(\mbf{k})  }{\omegab (\mbf{k})}  d \mbf{k}  .
\end{equation}
Here $ \omegab $ satisfies the following condition.

\begin{quote}
\textbf{(B.1)} $\omegab $ is continuous and non-negative.  
\end{quote}

$\;$ \\
As a remark, we consider  a physical example of $\omegab$. Let   
$ \omegab (\mbf{k}) = \sqrt{\mbf{k}^2 + m^2 }$, where $m \geq0 $ denotes the mass of the  field. 

 Let $ \ms{K_{\bos}} = \overline{\ms{K}_{\bos}^0}^{\| \cdot \|_{\ms{K}_{\bos} }} $.
From Proposition B,  there exists a Gaussian random process $ \{ \phi (f)  \}_{f \in \ms{K_{\bos}} } $   indexed by 
 $ \ms{K}_{\bos} $ on a probability space   $(Q_{\bos}, \mathfrak{B}_{\bos} , P_{\bos} )$. 
Let $ \check{\omegab} =  \omega_{\bos} (- i \nabla ) $.  The free Hamiltonian of the bose field is given by 
\[
H_{\bos} \; = \;    d \Gamma (  \check{\omegab}  ) .
\]
Then the  triplet  $ (   L^2 ( Q_{\ms{K}_{\bos}}  ) , \, \Hb,   
 \{\phi (f) \}_{ f \in  \ms{K}_{\bos} } )$ of the  scalar Bose field is constructed.
 
\subsection{Functional integral representation for scalar Bose fields}
In this subsection, we apply the Euclidean quantum filed theory. For the detail of this subject, refer to e.g. 
(\cite{Arai10}; Section 7) and (\cite{LHB11}; Section 5).
$\;$ \\ 
Let 
\[
\ms{K}_{\Euc}^0 \; = \; 
\left\{  f \in \ms{S}_{\real}' ( \mbf{R}^{1+d} ) \left| \frac{}{} \right. \int_{ \mbf{R}^{1+d}  }
\frac{|\hat{f} ( k_0, \mbf{k})|^2 }{\omegab (\mbf{k})^2 +  k_0^2 }  d k_0 d \mbf{k} \; < \;  \infty  \right\} ,
\]
 and   
\begin{equation}
( g , f )_{\ms{K}_{\Euc}} \; = \; 
 \int_{\Rd} 
\frac{ \overline{\hat{g} ( k_0 , \mbf{k})} \hat{f} ( k_0 , \mbf{k})  }{\omegab (\mbf{k})^2 +  k_0^2}  d k_0 d \mbf{k}  .
\end{equation}
Let $ \ms{K}_{\Euc} = \overline{\ms{K}^0_{\Euc}}^{\| \cdot \|_{ \ms{K}_{\Euc} }} $.
From Proposition B,  it is seen  that there exists  Gaussian random variables
 $ \{ \phi^{\Euc} (f) \}_{f \in \ms{K_{\Euc}} } $   indexed by  $ \ms{K}_{\Euc} $ 
on a probability space   $(Q_{\Euc}, \mathfrak{B}_{\Euc} , P_{\Euc} )$.

$\;$ \\
The relation between $\ms{K}_{\bos}$ and $\ms{K}_{\Euc}$ is as follows. 
For the  delta function $\delta_{t}  \in \ms{S}' (\mbf{R}) $ with $<\delta_t , \phi > = \phi (t) $,  it is seen that 
\begin{align}
& (g, e^{-|t-s | \check{\omegab}  }f )_{\ms{K}_{\bos}} \; = \; ( \delta_s \tens g , \delta_t \tens f   )_{\ms{K}_{\Euc}} 
, \qquad s \ne t  ,  \\
& \| f \|_{\ms{K}_{\bos}} \; = \; \| \delta_t  \tens  f  \|_{\ms{K}_{\Euc}} .
\end{align}
Then the  isometric operator $j_{t} : \ms{K}_{\bos}  \to \ms{K}_{\Euc}$ is defined by  
  $j_{t} f = \delta_{t}  \tens f  $.
  Let $J_{t} = \Gamma (j_t )$. Then it is seen that $e^{-t \Hb } = J_0^{\ast} J_{t}$ and 
\begin{align}
( \Phi , e^{-t \Hb }\Psi )_{L^2 (Q_{\bos})} \; = \; \mathbb{E}_{\Euc} [ (J_0 \Phi )^\ast (J_t \Psi ) ] ,
\end{align}
where $\mathbb{E}_{\Euc} (X) = \int_{Q_{\Euc}} X (\tilde{q}) d P_{\Euc} (\tilde{q}) $.
For $ D \subset \mbf{R}$, let us set $ E_{D} = \Gamma (e_{D})$ where $e_{D}$ is the projection  onto $  \ms{K}_{\Euc} ( D )   = \{  f \in  \ms{K}_{\Euc}  | supp f \in D \times \mbf{R}^d \}$.  It is seen that $E_{D} $ has the Markov property  such that    $ E_{[a, b]}E_{\{c \} } E_{[d, e]}  = E_{[a, b]} E_{[d, e]} $ for $a \leq  b \leq c \leq d \leq  e$.
Let $E_s  = J_s J_s^{\ast}$. Then it is known that $ E_{s} = E_{\{s \}}$. 
It is seen that 
\begin{equation}
  J_s G (\phi (f) ) J_s = E_s G (\phi^{\Euc} ({\delta_s \tens f}) ) E_s .  \label{8/10.3.1}
  \end{equation}  
  for $G \in   L^\infty ( \mbf{R}) $. 
Then   by using Trotter-Kato product formula, the following proposition holds 
(Refer to  e.g. \cite{Arai10}; Theorem 7.19).  
    
$\;$ \\   
\textbf{Proposition C} \\ 
Assume that $V_{\bos} $ is continuous function on $\Rd$ with bounded from bellow. 
Then
\[
( \Phi , e^{-t (H_{\bos} \dotplus  V_{\bos} (\phi (f) ) )} \Psi )_{L^2{ (Q}_{\bos}) } \; 
\;  = \; \; \mathbb{E}_{\Euc}
[  \overline{ (J_0 \Phi ) } (J_t \Psi ) e^{ - \int_0^t V_{\bos} (\phi^{\Euc} ( {\delta_s \tens f} ) )  ds }  ] . 
\]

\section{Main  Theorem and Proofs}
\subsection{Interacting system and main theorem}
The interaction system between the semi-relativistic particle and a scalar Bose fields is defined as follows. 
The state space for the system is given by
\[
\ms{H} \; = \; \LtwoRdx \tens L^2 (Q_{\bos } ) .
\]
 $\ms{H} $ can be decomposed as  $ \ms{H}  \simeq  \int^{\oplus}_{\mathbf{R}^d}  L^2 (Q_{\bos } ) d \mbf{x} $
 where $ \int^{\oplus}$  denotes the fibre direct integral. The total Hamiltonian of the system is defined by
 form sum of the free Hamiltonian and interaction 
\begin{equation}
H_{\kappa} \; \; = \; \; H_{0} \;  \dotplus \; \kappa \HI , \qquad \kappa \in \mbf{R} ,
\end{equation}
where $
H_{0} =  \Hpar  \tens I  +   I \tens H_{\bos} $ and
and  $\HI$ is given by  
\[
\HI \; = \; \int^{\oplus}_{\Rd} P (\phi ( {\rho_{\mbf{x}}}) ) d \mbf{x} 
\]
with $P(\lambda  ) = \sum\limits_{j=1}^{2n} c_{j} \lambda^j  $,
 $c_j \in \mbf{R}$, $j=1, \cdots , 2n-1$, $c_{2n} > 0$ and  $ \rho_{x}$ satisfying  the following conditions.  
\begin{quote}
\textbf{(A.1)} For each $\mbf{x} \in \Rd $, $f_{\mbf{x}} \in  \ms{K}_{\bos} $ and $
\sup_{\mbf{x} \in \Rd } \|  \rho_{\mbf{x}} \|_{\ms{K}_\bos} < \infty  $. \\
\textbf{(A.2)} For each $t \in \mbf{R}$, the map $ \mbf{R} \ni \mbf{x} \mapsto \delta_t \tens \rho_{\mbf{x}}
 \in \ms{K}_{\Euc} $  is strongly continuous. 
\end{quote}

$\;$ \\
For a  physical example of the interaction, let $ \rho_{\mbf{x} } (\mbf{y}) = \rho (\mbf{y} -\mbf{x}) $ for $\rho \in  \ms{S}'_{\textrm{real}} (\Rd )$. 
Then $ \hat{\rho}_{\mbf{x} } (\mbf{k}) = \hat{\rho} (\mbf{k}) e^{i \mbf{k} \cdot \mbf{x} }  $. Then we see that   the conditions 
 \textbf{(A.1)} and \textbf{(A.2)} are satisfied. The field operator 
$\phi (\rho_{\mbf{x}}) $  can be unrigorously represented as
\[
\phi (\rho_{\mbf{x}}) = \int_{\Rd} \frac{  \hat{\rho} (\mbf{k}) }{ \sqrt{ 2\omega (\mbf{k})}} 
\left( a_{\mbf{k}}  e^{-i \mbf{k} \cdot \mbf{x} } +  a_{\mbf{k}}^{\dagger} e^{i \mbf{k} \cdot \mbf{x} }  \right) d \mbf{k}
\]
where $a_{\mbf{k}} $ and $a^{\dagger}_{\mbf{k}}$ denote the kernel of an annihilation operator and creation operator, respectively.

$\;$ \\
We prepare for some notations. Let  $\ms{D}_0 \, = \, C_0^{\infty} (\Rd) \hat{\tens} \ms{D}_{0, \ms{K}_{\bos}}  $
 where $\hat{\tens} $ denotes the algebraic tensor product.
For $\Psi \in \ms{H} $,  we set $ \Psi_{\mbf{x}} (q)  = \Psi (\mbf{x} , q) $. 
Unless confusion arises, we identify $X \tens I$ with $X$ and $I \tens Y$ with $Y$. 
Let
\[
( \Omega_{\parti \times \Euc} , \mathfrak{B}_{\parti \times \Euc} , P^{\mbf{x}}_{\parti \times \Euc} ) = 
 ( \Omega_{\parti}  \times Q_{\Euc}  ,  \overline{\mathfrak{B}_{\parti} \times  
 \mathfrak{B}_{\Euc}} ,
  \overline{ P_{\parti}^{\mbf{x}} \times P_{\Euc}}  ) , 
\]
and we use the notation $\mathbb{E}^{\mbf{x}}_{\parti \times \Euc} [Z]  = \int_{\Omega_{\parti \times \Euc}} Z (\zeta  )d  P^{\mbf{x}}_{\parti \times \Euc} (\zeta)$.  

$\;$ \\
The main theorem in this paper is as follows.
\begin{theorem} 
 Assume \textbf{(S.1)},  \textbf{(B.1)}, \textbf{(A.1)} and \textbf{(A.2)}. Then it follows that
\[
(\Phi, e^{-tH_{\kappa} } \Psi )_\ms{H}
\; = \;
\int_{\Rd} \mathbb{E}_{\parti \times \Euc}^{\mbf{x}}
[  \overline{( J_0 \Phi_{\mbf{X}_0} )}  ( J_t \Psi_{\mbf{X}_t} ) e^{-\int_0^t V(\mbf{X}_s) ds } 
  e^{- \kappa   P (   \phi^{\Euc} (   \int_{0}^t  \delta_s  \tens \rho_{ {}_{\mbf{X}_s } }  ds) ) ) } ] d\mbf{x}.
\]
\end{theorem}
$\;$ \\
Now we consider an application of the above theorem. 
For  a self-adjoint  $H$ with bounded from below,  it is said that $ H $ has the ground state if  the infimum of the spectrum of $H$ is the eigenvalue.  It is seen that  $e^{-t \Hpar} $ and $e^{-t\Hb}$ are positivity improving operators. Hence from the above  functional integral representation, the  next corollary   immediately follows. 
\begin{corollary} $\;$  \\
Assume  \textbf{(S.1)}, \textbf{(B.1)}, \textbf{(A.1)} and \textbf{(A.2)}. Then  if $ H_{\kappa} $ has the ground state, 
 it is unique.  
\end{corollary}

\subsection{Proof of  Theorem 4.1}
To prove the Theorem 4.1, we show the following proposition.
\begin{proposition} 
Let $G_j \in L^{\infty} (\Rd) $, $j=1, \cdots n $. 
Then it follows that for $\Phi , \Psi  \in \ms{D}_0$, 
\begin{align*}
&(\Phi, e^{-t_1 H_{0} } G_1 (\phi ( {\rho_{\mbf{x}}} ) ) e^{-(t_2 - t_1 )H_{0} }  
 G_2 (\phi ({\rho_{\mbf{x}}} )  )
\cdots G_{n-1} (\phi ( {\rho_{\mbf{x}}}))
e^{-(t_{n} - t_{n-1} ) H_{0} }  \Psi ) \; \\
& = \; \int_{\Rd} \; \mathbb{E}_{\parti \times \Euc}^{\mbf{x}}
\left[  \overline{ (J_{0}  \Phi_{\mbf{X}_{t_1}} ) }  \left( \prod_{j=1}^{n-1}
G_j (\phi ( { \delta_{t_j } \tens \rho_{\mbf{x}}} ))   \right)(J_{t_n} \Psi_{\mbf{X}_{t_n} } ) 
e^{-\int_0^{t_n} V(\mbf{X}_s) ds } \right] d\mbf{x} . 
\end{align*}
\end{proposition}
\textbf{(Proof)} By using $e^{-(t-s)\Hb } = J_s^{\ast} J_{t}$ for $t>s$ and  (\ref{8/10.3.1}), it is seen that 
\begin{align}
& (\Phi, e^{-t_1 H_{0} } G_1 (\phi ( {\rho_{\mbf{x}}})) e^{-(t_2 - t_1 )H_{0} }  G_2 (\phi ( {\rho_{\mbf{x}}} ) ) 
\cdots G_{n-1} (\phi  (\rho_{\mbf{x}}))
e^{-(t_{n} - t_{n-1} ) H_{0} }  \Psi )_{\ms{H}}   \notag  \\
&= \int_{\Rd} (\Phi_{\mbf{x}} , e^{-t_1 H_{0} } G_1 (\phi ( {\rho_{\mbf{x}}})) e^{-(t_2 - t_1 )H_{0} }  
 G_2 (\phi ( {\rho_{\mbf{x}}} ))  \cdots G_{n-1} (\phi ( {\rho_{\mbf{x}}})) e^{-(t_{n} - t_{n-1} ) H_{0} }  \Psi_{\mbf{x}} )_{L^2 (Q_{\bos})} d \mbf{x}  \notag \\
&= \int_{\Rd} (\Phi_{\mbf{x}} , e^{-t_1 \Hpar }  
J_{0}^\ast \left( E_{t_{1}} G_1 (\phi^{\Euc} (  \delta_{t_1 } \tens \rho_{\mbf{x}} ))  E_{t_{1}} \right)
e^{-(t_2 - t_1 ) \Hpar } \left( E_{t_{2}}  G_2 (\phi^{\Euc} (  \delta_{t_2} \tens \rho_{\mbf{x}} ))  E_{t_{2}} \right)  
 \times  \notag    \\ 
& \qquad     \cdots   \times \left( E_{t_{n-1}}  G_{n-1} (\phi^{\Euc} ({ \delta_{t_{n-1}} \tens \rho_{\mbf{x}}} )  )E_{t_{n-1}} \right)    
 J_{t_n}  
 e^{-(t_{n} - t_{n-1} ) \Hpar }  \Psi_{\mbf{x}} )_{L^2 (Q_{\bos})} d \mbf{x}  . \label{8/7.1}
\end{align}
By using Markov property of $E_{t_j}$ and the Feynman-Kac formula (\ref{FKS_relativistic}), we have  
\begin{align*}
 (\ref{8/7.1}) &= \int_{\Rd}   \mathbb{E}_{\Euc} [ (\overline{ J_{0} \Phi_{\mbf{x}} }) e^{-t_1 \Hpar }  
  G_1 (\phi^{\Euc} ( { \delta_{t_1 } \tens \rho_{\mbf{x}}} ) )e^{-(t_2 - t_1 ) \Hpar }   
 G_2 (\phi^{\Euc} ( { \delta_{t_2} \tens  \rho_{\mbf{x}}} )) \times \\  
& \qquad \qquad  \cdots \times     G_{n-1} (\phi^{\Euc} ( { \delta_{t_{n-1}} \tens \rho_{\mbf{x}}} ))   e^{-(t_{n} - t_{n-1} ) \Hpar }  J_{t_n}   \Psi_{\mbf{x}} )  ] d \mbf{x} \\
&= \int_{\Rd}   \mathbb{E}_{\Euc} [ \mathbb{E}_{\parti}^{\mbf{x}} [ ( \overline{J_{0} \Phi_{\mbf{x}} })  
 \left( \prod_{j=1}^{n-1}  G_j (\phi^{\Euc} ({ \delta_{t_j } \tens \rho_{\mbf{X}_{t_j } }} ))  
\right)
   ( J _{t_n}   \Psi_{\mbf{X}_{t_n } } )   e^{-\int_{0}^{t_n} V (\mbf{X}_s )  ds } ] ] d\mbf{x}  . \\
   \end{align*}
Thus the proof is obtained. $\blacksquare$

$\;$ \\
\textbf{(Proof of Theorem 4.1)} \\
  Let $\Phi , \Psi  \in \ms{D}_{0}$. By Proposition  4.3 and Trotter-Kato product formula we have 
\begin{align}
(\Phi , e^{-t H_{\kappa } }\Psi ) 
&= \lim_{n \to \infty} (\Phi , (e^{-\frac{t}{n} H_0}  e^{ - \frac{t}{n} \kappa \HI })^n \Psi ) \notag  \\
&= \lim_{n \to \infty } \int_{\Rd} \mathbb{E}_{\parti \times \Euc}^{\mbf{x}} 
[   ( \overline{J_{0}   \Phi_{\mbf{X}_{0} } }) 
  ( J_{t}   \Psi_{\mbf{X}_{t} } )  
  e^{-\kappa 
P ( \phi^{\Euc} ( \sum\limits_{j=1}^n  \left( \frac{t}{n}  \right) {\delta_{tj/n} \tens \rho_{ \mbf{X}_{tj/n} } }  ) ) } e^{-\int_{0}^t V (\mbf{X}_s )  ds }  ] d\mbf{x} 
 \label{8/10.4.3} 
\end{align}
Here note that 
\begin{align}
  \| \delta_{t +\epsilon } \tens \rho_{\mbf{X}_{t + \epsilon}} 
  - \delta_{t} \tens \rho_{\mbf{X}_t}  \|_{\ms{K}_\Euc} 
&\leq 
 \|    (\delta_{t + \epsilon}  - \delta_t ) \tens \rho_{\mbf{X}_{t+\epsilon}} \|_{\ms{K}_{\Euc}}  + 
     \| \delta_{t } \tens ( \rho_{\mbf{X}_{t+\epsilon}}  -\rho_{\mbf{X}_t} ) \|_{\ms{K}_{\Euc}}  \notag \\
& \leq     \|    ( 1- e^{\epsilon \check{\omegab}} )^{1/2} \rho_{\mbf{X}_{t+\epsilon}} \|_{\ms{K}_{\bos}} + 
   \|   (\rho_{ \mbf{X}_{t+ \epsilon}}  - \rho_{ \mbf{X}_t} ) \|_{ \ms{K}_{\bos}} \label{8/10.6}
 \end{align}
It is known that the map 
 $ s \to \mbf{X}_s (\xi ) $ is continuous  for each $\xi  \in \Omega_{\parti}$ except finite points.  
Then from \textbf{(A.2)} and (\ref{8/10.6}), the map $\mbf{R} \ni t  \mapsto \delta_{t} \tens \rho_{\mbf{X}_t}  \in \ms{K}_{\Euc}$ is strongly  continuous almost surely.   Then from this continuity and (\ref{8/10.4.3}),  we have
 \begin{equation} 
(\Phi , e^{-t H_{\kappa } }\Psi )  = \int_{\Rd} \mathbb{E}_{\parti \times \Euc}^{\mbf{x}}
[  ( \overline{J_0 \Phi_{\mbf{X}_0} })  ( J_t \Psi_{\mbf{X}_t} )
   e^{- \kappa   P (   \phi^{\Euc} (   \int_{0}^t  \delta_s  \tens \rho_{ {}_{\mbf{X}_s } }  ds) ) ) }
  e^{-\int_0^t V(\mbf{X}_s) ds }  ] d \mbf{x} . 
\end{equation}
Since $\ms{D}_0$ is  dense in $\ms{H}$,  the proof is obtained. $\blacksquare$

$\quad$ \\
{\Large Acknowledgments} \\
It is a pleasure to thank Professor  Fumio Hiroshima for his advice and comments.

\end{document}